\documentclass[11pt]{article}
\usepackage{enumerate}
\usepackage{amssymb,a4wide,latexsym,makeidx,epsfig,fleqn}
\usepackage{amsthm}
\usepackage{amsmath}
\usepackage{enumerate}
\usepackage{graphicx}
\usepackage{float}
\usepackage{indentfirst}
\usepackage[numbers,sort&compress]{natbib}
\usepackage{color}
\allowdisplaybreaks[4]
\newtheorem{theorem}{Theorem}[section]

\newtheorem{lemma}{Lemma}[section]

\newtheorem{problem}{Problem}[section]

\begin{document}
\textwidth 150mm \textheight 225mm
\title{On the $\alpha$-index of minimally $k$-(edge-)connected graphs for small $k$ \thanks{Supported by the National Natural Science Foundation of China (No. 12271439).}}
\author{{Jiayu Lou$^{a,b}$, Ligong Wang$^{a,b,}$\thanks{Corresponding author.\newline
\setlength{\hangindent}{1.9em}
E-mail addresses: lgwangmath@163.com, jyloumath@163.com, ym19980508@mail.nwpu.edu.cn.}, Ming Yuan$^{a}$}\\
{\small $^a$School of Mathematics and Statistics, Northwestern
Polytechnical University,}\\ {\small  Xi'an, Shaanxi 710129,
P.R. China.}\\
{\small $^b$ Xi'an-Budapest Joint Research Center for Combinatorics, Northwestern
Polytechnical University,}\\
{\small Xi'an, Shaanxi 710129,
P.R. China. }\\}
%{\small E-mail address: jyloumath@163.com, lgwangmath@163.com, ym19980508@mail.nwpu.edu.cn} }
\date{}
\maketitle
\begin{center}
\begin{minipage}{120mm}
\vskip 0.3cm
\begin{center}
{\small {\bf Abstract}}
\end{center}
{\small Let $G$ be a graph with adjacency matrix $A(G)$ and let $D(G)$ be the diagonal matrix of vertex degrees of $G$. For any real $\alpha \in [0,1]$, Nikiforov defined the $A_\alpha$-matrix of a graph $G$ as $A_\alpha(G)=\alpha D(G)+(1-\alpha)A(G)$.
The largest eigenvalue of $A_\alpha(G)$ is called the $\alpha$-index or the $A_\alpha$-spectral radius of $G$. A graph is minimally $k$-(edge)-connected if it is $k$-(edge)-connected and deleting any arbitrary chosen edge always leaves a graph which is not $k$-(edge)-connected.
In this paper, we characterize the minimally 2-edge-connected graphs and minimally 3-connected graph with given order having the maximum $\alpha$-index for $\alpha \in [\frac{1}{2},1)$, respectively.

 }
\end{minipage}
\end{center}

%%%%%%%%%%%%%%%%%%%%%%%%%%%%%%%%%%%%%%%%%%%%%%%%%%%%%%%%%%%%%%%%%%%%%%%%%%%%%%
%%% Part 1: Introduction %%%
\section{Introduction }
\label{sec:ch6-introduction}
% Para1: basic concept
\noindent Throughout this paper, we only consider simple, finite and undirected graphs. Let $G=(V(G),E(G))$ be a graph with vertex set $V(G)=\{v_1,v_2,\ldots,v_n\}$ and edge set $E(G)$. The order of $G$ is the number $n=|V(G)|$ of its vertices and the size of $G$ is the number $m=|E(G)|$ of its edges.
The neighborhood $N_G(v)$ (or, $N(v)$ for short) of a vertex $v$ is $\{u\in V(G):uv\in E(G)\}$ and the number $d_G(v)=|N_G(v)|$ (or, $d(v)$ for short) is the degree of $v$. The closed neighbor set is defined as $N_G[v]=N_G(v)\cup \{v\}$ (or, $N[v]$ for short).
Let $\Delta(G)$ and $\delta(G)$ be the maximum degree and minimum degree of $G$, respectively. The average degree of the neighbors of a vertex $u$ of $G$ is $m(u)=\frac{\sum_{u v \in E(G)} d(v)}{d(u)}$.
For graphs $F$ and $G$, we use $H\subseteq G$ to denote that $G$ contains $H$ as a subgraph.
For a vertex set $S \subseteq V(G)$, let $G[S]$ be the subgraph of $G$ induced by $S$. For $A, B \subset V(G)$, we denote by $e(A)$ the number of edges in $G[A]$ and by $e(A, B)$ the number of edges with one endpoint in $A$ and one endpoint in $B$.
A cut vertex (resp., cut edge) of a graph is a vertex (resp., edge) whose deletion increases the number of components of a graph. A cycle $C$ of $G$ is said to have a chord if there is an edge of $G$ that joins a pair of non-adjacent vertices from $C$.

Denote by $P_n$, $C_n$ and $K_{a,n-a}$ the path, the cycle and the complete bipartite graph on $n$ vertices, respectively. For two graphs $G_1$ and $G_2$, we define $G_1\cup G_2$ to be their disjoint union. The join $G_1\vee G_2$ is obtained from $G_1\cup G_2$ by joining every vertex of $G_1$ with every vertex of $G_2$ by an edge. Let $W_n=K_1\vee C_{n-1}$ be the wheel graph on $n$ vertices.
For $k\geq1$, denote by $F_k$ the graph on $2k+1$ vertices consisting of $k$ triangles which intersect in exactly one common vertex.

The adjacency matrix of $G$ is defined as the matrix $A(G)=(a_{u,v})_{u,v\in V(G)}$ with $a_{u,v}=1$ if $u$, $v$ are adjacent, and $a_{u,v}=0$ otherwise. Let $D(G)$ be the diagonal matrix of vertex degrees of $G$. The Laplacian matrix and signless Laplacian matrix of $G$ are defined as $L(G)=D(G)-A(G)$ and $Q(G)=D(G)+A(G)$, respectively. The largest eigenvalue of $A(G)$ is called the index or the spectral radius of $G$, and the largest eigenvalue of $Q(G)$ is called the $Q$-index or the signless Laplacian spectral radius of $G$. The largest eigenvalue of $L(G)$ is called the $L$-index or the Laplacian spectral radius of G.
To track the gradual change of $A(G)$ into $Q(G)$, Nikiforov \cite{N2017} proposed to study the convex linear combination $A_\alpha(G)$ of $A(G)$ and $D(G)$ defined by
$$A_\alpha(G)=\alpha D(G)+(1-\alpha)A(G)$$ for any real $\alpha \in [0,1]$.
Note that $A(G)=A_0(G)$, $D(G)=A_1(G)$ and $Q(G)=2A_{\frac{1}{2}}(G).$ The largest eigenvalue of $A_\alpha(G)$, denoted by $\rho_\alpha(G)$, is called the $\alpha$-index or the $A_\alpha$-spectral radius of $G$. Denote by $\lambda(M)$ the largest eigenvalue of a matrix $M$. For a connected graph $G$, $A_\alpha(G)$ is irreducible.
By the Perron-Frobenius Theorem, $\rho_\alpha(G)$ is positive, and there exists a unique positive unit eigenvector corresponding to $\rho_\alpha(G)$, which is called the $\alpha$-Perron vector of $G$.

Let $\mathcal{G}$ be a specified class of graphs. Brualdi and Solheid \cite{BS} proposed the problem of characterizing the graphs for which the maximal or the minimal spectral invariants are attained in $\mathcal{G}$. For
many different $\mathcal{G}$ and different spectral invariants (such as index, L-index and Q-index), this problem has been extensively studied (see \cite{DS,LC,LGW,Z}). For $A_\alpha$-index, there are some relative results in \cite{LYL,CH,LCG}. We refer the interested reader to the surveys \cite{CZ,LLF,N2011} for more results.

A graph is $k$-connected (resp., $k$-edge-connected) if removing fewer than $k$ vertices (resp., edges) always leaves the remaining graph connected, and is minimally $k$-(edge-)connected if it is $k$-connected (resp., $k$-edge-connected) and deleting any arbitrary chosen edge always leaves a graph which is not $k$-connected (resp., $k$-edge-connected).
In recent works, some researchers restrict $\mathcal{G}$ to (minimally) $k$-(edge-)connected graphs of order $n$ or size $m$. A graph is minimally 1-(edge-)connected if and only if it is a tree. It is natural to consider the problem of finding the graph with the maximum ($Q$-)index among all minimally $k$-(edge-)connected graphs for $k\geq 2$. Fan, Goryainov and Lin \cite{FGL} asked the following question for $k\geq 2$.
\noindent\begin{problem}
What is the maximum ($Q$-)index and what are the corresponding extremal graphs among minimally $k$-(edge-)connected graph for $k\geq 2$?
\end{problem}

Chen and Guo \cite{CG} characterized the minimally 2-(edge-)connected graphs with given order having the maximum index. Lou, Min and Huang \cite{LMH} characterized the minimally 2-(edge-)connected graphs with given size having the maximum index. Fan, Goryainov and Lin \cite{FGL} not only determined the minimally 2-(edge-)connected graphs with given order having the maximum $Q$-index, but characterized the minimally 3-connected graphs with given order having the maximum ($Q$-)index. Guo and Zhang \cite{GZ2} characterized the (minimally) 2-connected graphs with given size having the maximum $Q$-index. Zhang and Guo \cite{ZG} determined the (minimally) 2-connected graphs with given size having the maximum $L$-index. Lou, Wang and Yuan \cite{me} determined the minimally 2-connected graphs with given order or size having the maximum $\alpha$-index, respectively. They also asked the following question for $\alpha$-index.
% Problem
\noindent\begin{problem}
What is the maximum $\alpha$-index and what are the corresponding extremal graphs among minimally $k$-(edge)-connected graph for $k\geq 2$?
\end{problem}

In this paper, we characterize the minimally 2-edge-connected graphs and minimally 3-connected graph with given order having the maximum $\alpha$-index for $\alpha \in [\frac{1}{2},1)$, respectively.
\noindent\begin{theorem}\label{th1.1}
Let $G$ be a minimally 2-edge-connected graph with order $n$ and $\alpha \in [\frac{1}{2},1)$.
\begin{enumerate}[{\rm(i)}]
\item If $n \geq 7$ is an odd number, then $\rho_{\alpha}(G)\leq \rho_{\alpha}(F_{\frac{n-1}{2}})$, with equality if and only if $G\cong F_{\frac{n-1}{2}}$.
\item If $n \geq 8$ is an even number, then $\rho_{\alpha}(G)\leq\rho_{\alpha}(K_{2,n-2})$, with equality if and only if $G\cong K_{2,n-2}$.
\end{enumerate}
\end{theorem}
\noindent\begin{theorem}\label{th1.2}
Let $G$ be a minimally 3-connected graph with order $n\geq7$. If $\alpha \in [\frac{1}{2},1)$,
%\begin{enumerate}[{\rm(i)}]
%\item If $n \geq 10$ is an even number, then $\rho_{\alpha}(G)\leq\rho_{\alpha}(K_{2,n-2})$, with equality if and only if $G\cong K_{2,n-2}$.
then $\rho_{\alpha}(G)\leq \rho_{\alpha}(W_n)$, with equality if and only if $G\cong W_n$.
%\end{enumerate}
\end{theorem}
\par{The remainder of the paper is organized as follows. In Section 2, we recall some useful lemmas that will be used later. In Sections 3 and 4, we give proofs of Theorems \ref{th1.1} and \ref{th1.2}, respectively.}
%%% Part 1 end %%%
%%%%%%%%%%%%%%%%%%%%%%%%%%%%%%%%%%%%%%%%%%%%%%%%%%%%%%%%%%%%%%%%%%%%%%%%%%%%%%

%%%%%%%%%%%%%%%%%%%%%%%%%%%%%%%%%%%%%%%%%%%%%%%%%%%%%%%%%%%%%%%%%%%%%%%%%%%%%%
%%% Part 2: Preliminaries %%%
\section{Preliminaries}
\label{sec:Preliminaries}
\noindent In this section, we introduce some useful lemmas that are used in the proof of our main results.

\noindent\begin{lemma}\label{ub1} {\normalfont(\cite{N2017})}
If $G$ is a graph with no isolated vertices, then
$$\rho_{\alpha}(G) \leq \max _{u \in V(G)}\left\{\alpha d(u)+\frac{1-\alpha}{d(u)} \sum_{u v \in E(G)} d(v)\right\}.$$
If $\alpha \in (\frac{1}{2},1)$ and $G$ is connected, with equality if and only if $G$ is regular.
\end{lemma}

\noindent\begin{lemma}\label{ub2} {\normalfont(\cite{WWT})}
For any graph $G$ and $\alpha\in[0, 1)$. Then
$$
\rho_{\alpha}(G) \leq \max _{uv \in E(G)}\left\{\frac{\alpha\left(d(u)+d(v)\right)+\sqrt{\alpha^2\left(d(u)-d(v)\right)^2+4(1-\alpha)^2 m(u) m(v)}}{2}\right\} .
$$
If $G$ is connected, the equality holds if and only if $G$ is regular or bipartite semi-regular.
\end{lemma}

\noindent\begin{lemma}\label{lb} {\normalfont(\cite{N2017})}
Let G be a graph with $\Delta(G)=\Delta$. If $\alpha \in[0,\frac{1}{2}]$, then
$$\rho_{\alpha}(G) \geq \alpha(\Delta+1).$$
If $\alpha \in[\frac{1}{2},1)$, then
$$\rho_{\alpha}(G) \geq \alpha \Delta+\frac{(1-\alpha)^{2}}{\alpha}.$$
\end{lemma}

\noindent\begin{lemma}\label{degree} {\normalfont(\cite{BM})}
If $G$ is a minimally $k$-(edge-)connected graph, then $\delta(G)=k.$
\end{lemma}

\noindent\begin{lemma}\label{mn} {\normalfont(\cite{BM})}
If $G$ is a minimally 2-edge-connected graph with order $n$, then $m\leq2n-2$.
\end{lemma}

\noindent\begin{lemma}\label{chord} {\normalfont(\cite{LMH})}
If $G$ is a minimally 2-edge-connected graph, then no cycle of $G$ has a chord.
\end{lemma}

%\noindent\begin{lemma}\label{3m} {\normalfont(\cite{B})}
%If $G$ is a minimally 3-connected graph with order $n$, then $m\leq3n-9$, with equality if and only if $G\cong K_{3,n-3}$.
%\end{lemma}

\noindent\begin{lemma}\label{degree3} {\normalfont(\cite{B})}
If $G$ is a minimally 3-connected graph with order $n$, then each cycle contains at least two vertices of degree 3.
\end{lemma}

\noindent\begin{lemma}\label{switch} {\normalfont(\cite{NR,XLL})}
Let $G$ be a connected graph with $\alpha \in[0,1)$. For $u, v \in V(G)$, suppose $N \subseteq N(v) \backslash(N(u) \cup\{u\})$. Let $G^{\prime}=G-\{v w: w \in N\}+\{u w: w \in N\}$. Let $X=\left(x_1, x_2, \ldots, x_n\right)^{T}$ be the $\alpha$-Perron vector of $G$ corresponding to $\rho_\alpha(G)$. If $N \neq \emptyset$ and $x_u \geq x_v$, then $\rho_\alpha\left(G^{\prime}\right)>\rho_\alpha(G)$.
\end{lemma}

\noindent\begin{lemma}\label{M} {\normalfont(\cite{N2017})}
Let $G_1$ be a $r_1$-regular graph of order $n_1$, and $G_2$ be a $r_2$-regular graph of order $n_2$. Then the largest eigenvalue of $A_\alpha(G_1\vee G_2)$ is
$$
\rho_\alpha(G_1\vee G_2)=\lambda\left(\begin{array}{cc}
r_1+\alpha n_2 & (1-\alpha)^2 n_1 n_2 \\
1 & r_2+\alpha n_1
\end{array}\right).
$$
\end{lemma}

%\noindent\begin{lemma}\label{lem7} {\normalfont(\cite{NR,XLL})}
%Let $G$ be a connected graph with $\alpha \in[0,1)$. For $u, v \in V(G)$, suppose $N \subseteq N(v) \backslash(N(u) \cup\{u\})$. Let $G^{\prime}=G-\{v w: w \in N\}+\{u w: w \in N\}$. Let $X=\left(x_1, x_2, \ldots, x_n\right)^{T}$ be the $\alpha$-Perron vector of $G$ corresponding to $\rho_\alpha(G)$. If $N \neq \emptyset$ and $x_u \geq x_v$, then $\rho_\alpha\left(G^{\prime}\right)>\rho_\alpha(G)$.
%\end{lemma}

%We say that $u$ and $v$ are equivalent in $G$ if there exists an automorphism $p: G\rightarrow G$ such that $p(u)=v$.
%\noindent\begin{lemma}\label{lem8} {\normalfont(\cite{N2017})}
%Let $G$ be a connected graph of order $n$, and let $u$ and $v$ be equivalent vertices in G. If $X=\left(x_{1}, x_{2}, \ldots, x_{n}\right)^{T}$ is an eigenvector to $\rho_{\alpha}(G)$, then $x_{u}=x_{v}$.
%\end{lemma}

\noindent\begin{lemma}\label{kab}
{\normalfont(\cite{N2017})}
Let $a \geq b \geq 1$. If $\alpha \in[0,1]$, then the largest eigenvalue of $A_\alpha\left(K_{a, b}\right)$ is
$$
\rho_\alpha(K_{a,b})=\frac{1}{2}(\alpha(a+b)+\sqrt{\alpha^2(a+b)^2+4 a b(1-2 \alpha)}).
$$
\end{lemma}

It is obvious that $\rho_\alpha(K_{2,n-2})=
\frac{1}{2}(\alpha n+\sqrt{\alpha^2 n^2+8(1-2\alpha)(n-2)})$.

\noindent\begin{lemma}\label{f}
Let $n$ be an odd integer. If $\alpha \in[0,1]$, then the largest eigenvalue of $A_\alpha(F_{\frac{n-1}{2}})$ is
$$
\rho_\alpha(F_{\frac{n-1}{2}})=\frac{1}{2}(\alpha n+1+\sqrt{\alpha^2n^2-10\alpha n+12\alpha+4n-3}).
$$
\end{lemma}
\noindent\textbf{Proof.}
Since $F_{\frac{n-1}{2}}=K_1 \vee \frac{n-1}{2}K_2$, then by Lemma \ref{M}, we have $$\rho_\alpha(F_{\frac{n-1}{2}})=\lambda\left(\begin{array}{cc}
\alpha(n-1) & (1-\alpha)^2(n-1) \\
1           &  1+\alpha
\end{array}\right).$$
By computation, we have
$$
\rho_\alpha(F_{\frac{n-1}{2}})=\frac{1}{2}(\alpha n+1+\sqrt{\alpha^2n^2-10\alpha n+12\alpha+4n-3}).\quad \qedsymbol
$$

%\noindent\begin{lemma}\label{W}
%If $\alpha \in[0,1]$, then the largest eigenvalue of $W_n$ is
%$$
%\rho_\alpha(W_n)=\frac{1}{2}(\alpha n+2+\sqrt{(\alpha n+2)-4(n-1)(4\alpha-1)}).
%$$
%\end{lemma}
%\noindent Proof. Since $W_n=K_1 \vee C_{n-1}$, then by Lemma \ref{M}, we have $$\rho_\alpha(W_n)=\lambda\left(\begin{array}{cc}
%\alpha(n-1) & (1-\alpha)^2(n-1) \\
%1           &  2+\alpha
%\end{array}\right).$$
%By computation, we have
%$$
%\rho_\alpha(W_n)=\frac{1}{2}(\alpha n+2+\sqrt{(\alpha n+2)-4(n-1)(4\alpha-1)}).\quad \qedsymbol
%$$

%\noindent\begin{lemma}\label{kf}
%$$
%\rho_\alpha(K_{2,n-2})< \rho_\alpha(F_{\frac{n-1}{2}})
%$$
%for $n>3$.
%\end{lemma}
%\noindent{\color{blue}\textbf{Proof.}}
%%% Part 2 end %%%
%%%%%%%%%%%%%%%%%%%%%%%%%%%%%%%%%%%%%%%%%%%%%%%%%%%%%%%%%%%%%%%%%%%%%%%%%%%%%%

%%%%%%%%%%%%%%%%%%%%%%%%%%%%%%%%%%%%%%%%%%%%%%%%%%%%%%%%%%%%%%%%%%%%%%%%%%%%%%
%%% Part 3: Proof Theorem 1.1 %%%
\section{Proof of Theorem \ref{th1.1}}
\label{sec:ch-sufficient}
\noindent In this section, we give the proof of Theorem \ref{th1.1}.\\ \hspace*{\fill} \\
\textbf{Proof.} \textbf{(\romannumeral1)}
Let $G$ be the graph having the maximum $\alpha$-index among all minimally 2-edge-connected graphs with odd order $n\geq7$. Then $\Delta(G)\leq n-1$. Notice that $F_{\frac{n-1}{2}}$ is a minimally 2-edge-connected graph. It follows that $\rho_\alpha(G)\geq \rho_\alpha(F_{\frac{n-1}{2}})$. By Lemma \ref{lb}, we have
\begin{equation}{\label{eq1}}
\rho_{\alpha}(F_{\frac{n-1}{2}})> \alpha\Delta(F_{\frac{n-1}{2}})+\frac{(1-\alpha)^{2}}{\alpha}
=\alpha(n-1)+\frac{(1-\alpha)^{2}}{\alpha}
\end{equation}
for $\alpha\in[\frac{1}{2},1)$. \\
Let $w$ be a vertex of $G$ such that
$$\alpha d(w)+\frac{1-\alpha}{d(w)} \sum_{wv \in E(G)}d(v)=\max_{u \in V(G)}\left\{\alpha d(u)+\frac{1-\alpha}{d(u)} \sum_{uv \in E(G)} d(v)\right\}.$$
By Lemma \ref{ub1}, we have
\begin{equation}{\label{eq2}}
\rho_{\alpha}(G)\leq \alpha d(w)+\frac{1-\alpha}{d(w)} \sum_{wv \in E(G)}d(v).
\end{equation}
Let $A_1$ be the isolated vertex set of $G[N(w)]$. Denote by $B=V(G)\backslash N[w]$ and $A_2=N(w)\backslash A_1$. Now we give Claims 1-3 to complete the proof of Theorem \ref{th1.1} (i).\\
\textbf{Claim 1.} $G[N(w)]=sK_2\cup tK_1$, where $s,t\geq 0$.\\
\textbf{Proof.}
Firstly, we assert that $G[N(w)]$ is $P_3$-free. Otherwise, there exists a $P_3=v_1v_2v_3$ in $G[N(w)]$ such that $wv_2$ is a chord of cycle $wv_1v_2v_3w$, which contradicts with Lemma \ref{chord}.
Secondly, we assert that $G[N(w)]$ is $C_l$-free, where $l\geq 3$. Otherwise, there exists a wheel graph $W_{l+1}$ in $G$ such that there is a cycle containing a chord, a contradiction.
Hence, we have $G[N(w)]=sK_2\cup tK_1$, where $s, t\geq 0$. $\qedsymbol$\\
\textbf{Claim 2.} $N_B(u)=\varnothing$ for any $u \in A_2$.\\
\textbf{Proof.}
Otherwise, there exists a vertex $v \in B$ that is adjacent to a vertex $u_2 \in A_2$. We assume that $u'_{2}\in N_{A_2}(u_2)$. If $N_B(v)=\varnothing$, then there exists a vertex $u \in A$ adjacent to $v$ by Lemma \ref{degree}. We consider three cases regarding vertex $u$. If $u \in A_1$, then a chord $wu_2$ exists in cycle $wu'_{2}u_2vuw$, a contradiction. If $u \in A_2$ and $u=u'_{2}$, then $u'_{2}u_2$ is a chord of the cycle $wu'_{2}vu_2w$, a contradiction. If $u\in A_2$ and $u\neq u'_{2}$, then there exists a chord $wu_2$ in cycle $wu'_{2}u_2vuw$, a contradiction. Hence we have $N_B(v)\neq\varnothing$.

Denote by $P(u\rightarrow v)$ the path from $u$ to $v$, then there exists a path $P(v\rightarrow v_t)$ in $G[B]$ such that $v_t$ is adjacent to a vertex $u'\in A$. Otherwise $u_2v$ will be a cut edge. If $u' \in A_1$, then a chord $wu_2$ exists in cycle $wu'_{2}u_2P(v\rightarrow v_t)u'w$, a contradiction. If $u' \in A_2$ and $u'=u'_{2}$, then $u'_{2}u_2$ is a chord of cycle $wu_{2}P(v\rightarrow v_t)u'w$, a contradiction. If $u'\in A_2$ and $u'\neq u'_{2}$, then $G'=G-u_2v-u_2v_t+wv+wv_t$ is a minimally 2-edge-connected graph. By Lemma \ref{switch}, we have $\rho(G')>\rho(G)$, which contradicts the maximality of $G$.
$\qedsymbol$\\
\textbf{Claim 3.} $\Delta(G)=n-1$.\\
\textbf{Proof.}
If $\Delta(G)\leq n-2$, then $|B|\geq1$. By Lemma \ref{degree}, we have $2\leq d(w)\leq \Delta(G)\leq n-2$. We consider the following two cases.

\textbf{Case 1.} $e(N(w))=0$.
In this case, we have $e(N(w),B)\leq d(w)|B|$. It follows that
$$
\begin{aligned}
\sum_{wv \in E(G)} d(v)
&=d(w) + 2e(N(w))+ e(N(w),B)
\leq d(w)(|B|+1).
\end{aligned}
$$
Combining this with (\ref{eq2}), we have
$$
\begin{aligned}
\rho_\alpha(G) &\leq \alpha d(w)+(1-\alpha)(|B|+1)=\alpha d(w)+(1-\alpha)(n-d(w))\\
&=\alpha d(w)+(1-\alpha)n=(2\alpha-1)d(w)+(1-\alpha)n\\
&\leq (2\alpha-1)(n-2)+(1-\alpha)n.
\end{aligned}
$$
Noting that
\begin{equation}{\label{eq3}}
\alpha(n-1)+\frac{(1-\alpha)^2}{\alpha}-((2\alpha-1)(n-2)+(1-\alpha)n)
=\frac{(2\alpha-1)^2}{\alpha}\geq0
\end{equation}
for $n\geq7$ and $\alpha \in[\frac{1}{2},1)$, we have
$\rho_\alpha(G)\leq\alpha(n-1)+\frac{(1-\alpha)^2}{\alpha}.$
Combining this with (\ref{eq1}), we have
$\rho_\alpha(G)< \rho_\alpha(F_{\frac{n-1}{2}})$, a contradiction.

\textbf{Case 2.} $e(N(w))\geq1$.

By Claims 1 and 2, we have $e(N(w),B)\leq(d(w)-2e(N(w)))|B|$. It follows that
$
\sum_{wv \in E(G)} d(v) \leq d(w)+2e(N(w))+(d(w)-2e(N(w)))|B|.
$
Combining this with (\ref{eq2}) and $|B|\geq1$, we have
$$
\begin{aligned}
\rho_\alpha(G) &\leq \alpha d(w)+(1-\alpha)\left(d(w)+2e(N(w))+(d(w)-2e(N(w)))|B|\right)\\
&=\alpha d(w)+(1-\alpha)(|B|+1)+\frac{2(1-\alpha)e(N(w))}{d(w)}(1-|B|)\\
&\leq\alpha d(w)+(1-\alpha)(n-d(w))\\
&\leq (2\alpha-1)(n-2)+(1-\alpha)n.
\end{aligned}
$$
Combining with (\ref{eq1}) and (\ref{eq3}), we have
$\rho_\alpha(G)< \rho_\alpha(F_{\frac{n-1}{2}})$, a contradiction.

Combining the above arguments, we have $\Delta(G)=n-1$. $\qedsymbol$

By Claim 3, we have $\Delta(G)=1$. Notice that $G[N(w)]$ consists of some edges and isolated vertices. By the maximality of $\rho_\alpha(G)$, we have $G\cong F_\frac{n-1}{2}.$

\textbf{(\romannumeral2)}  Let $G$ be a minimally 2-edge-connected graph with even order $n\geq8$. Then $\Delta(G)\leq n-2$ by Lemmas \ref{degree} and \ref{chord}. Notice that $K_{2,n-2}$ is a minimally 2-edge-connected graph. %By Lemma \ref{kab}, we have
%$$
%\begin{aligned}
%\rho_{\alpha}(K_{2,n-2})=\frac{1}{2}\left(\alpha n+\sqrt{\alpha^2 n^2+8(n-2)(1-2\alpha)}\right).
%\end{aligned}
%$$
By Lemma \ref{lb}, we have
\begin{equation}{\label{eq4}}
\rho_{\alpha}(K_{2,n-2})> \alpha\Delta(K_{2,n-2})+\frac{(1-\alpha)^{2}}{\alpha}
=\alpha(n-2)+\frac{(1-\alpha)^{2}}{\alpha}
\end{equation}
for $\alpha\in[\frac{1}{2},1)$.

If $\Delta(G)=n-2$, then $|B|=1$. If $e(N(w))=0$, then $G\cong K_{2,n-2}$ by Lemma \ref{degree}. If $e(N(w))\geq 1$, we assume that $B=\{v\}$ and $u_i\in A_i$ for $i = 1,2$. By analysis, when $e(N(w))=1$, the vertices $w$ and $u_2$ make the right-hand side of the inequality in Lemma \ref{ub2} obtain the maximum value. It follows that
$$
\begin{aligned}
\rho_{\alpha}(G) &\leq \max _{uv \in E(G)}\left\{\frac{\alpha\left(d(u)+d(v)\right)+\sqrt{\alpha^2\left(d(u)-d(v)\right)^2+4(1-\alpha)^2 m(u) m(v)}}{2}\right\}\\
&=\frac{\alpha\left(d(w)+d(u_2)\right)+\sqrt{\alpha^2\left(d(w)-d(u_2)\right)^2+4(1-\alpha)^2 m(w) m(u_2)}}{2}\\
&=\frac{\alpha n+\sqrt{\alpha^2(n-4)^2+8(1-\alpha)^2(n-3)}}{2}.
\end{aligned}
$$ for $n\geq8$ and $\alpha \in[\frac{1}{2},1)$.

Noting that $\alpha^2(n-4)^2+8(1-\alpha)^2(n-3)-(\alpha^2 n^2+8(1-2\alpha)(n-2))=-8(\alpha-1)^2<0$, we have
$$
\frac{\alpha n+\sqrt{\alpha^2(n-4)^2+8(1-\alpha)^2(n-3)}}{2}
<\frac{\alpha n+\sqrt{\alpha^2 n^2+8(1-2\alpha)(n-2)}}{2}
$$ for $n\geq8$ and $\alpha \in[\frac{1}{2},1)$. By Lemma \ref{kab}, we have $\rho_{\alpha}(G)<\rho_\alpha(K_{2,n-2})$.

Next, we consider $\Delta(G)\leq n-3$. By Lemma \ref{degree}, we have $2\leq d(w) \leq \Delta(G)\leq n-3$.
%By Lemma \ref{lem2}, we have
%\begin{equation}{\label{1}}
%\rho_{\alpha}(K_{2,n-2})\geq \alpha\Delta(K_{2,n-2})+\frac{(1-\alpha)^{2}}{\alpha}
%=\alpha(n-2)+\frac{(1-\alpha)^{2}}{\alpha}
%\end{equation}
%for $\alpha\in[\frac{1}{2},1)$. Thus we assume that $\Delta(G)\leq n-3$.
\\
By Lemma \ref{ub1}, we have
\begin{equation}{\label{eq5}}
\rho_{\alpha}(G)\leq \alpha d(w)+\frac{1-\alpha}{d(w)} \sum_{wv \in E(G)}d(v).
\end{equation}
%\begin{equation}{\label{3}}
%\begin{aligned}
 %\sum_{wv \in E(G)} d(v)=2 e(N(w))+e(N(w), V(G) \backslash N(w)).
 %\end{aligned}
%\end{equation}

Now, we prove that $\rho_\alpha(G)< \rho_\alpha(K_{2,n-2})$ for $2\leq d(w)\leq n-3$. We consider the following two cases.
%%case1

\textbf{Case 1.} $d(w)=2$.

\textbf{Subcase 1.1.} $e(N(w))=0$.

If $e(B)=0$, then $G\cong K_{2,n-2}$ by Lemma \ref{degree}. If $e(B)\geq1$, then there exists at least one path $P_n$ ($n\geq2$) in $B$. Since every cycle in $G$ has no chord, only the endpoints in $P_n$ are connected to each of the two vertices in $N(w)$,
and the remaining vertices are not connected to any vertex of $N(w)$.
It follows that $e(N(w),B)\leq2+d(w)(|B|-2)=2n-8$. Thus, $\sum_{wv \in E(G)} d(v)=d(w) + 2e(N(w)) + e(N(w),B)\leq 2n-6$. Combining this with (\ref{eq5}), we have
$$
\rho_\alpha(G) \leq 2 \alpha+\frac{1-\alpha}{2} \sum_{w v \in E(G)} d(v) \leq 2 \alpha+(1-\alpha)(n-3).
$$
Noting that
$$
\alpha(n-2)+\frac{(1-\alpha)^2}{\alpha}-(2 \alpha+(1-\alpha)(n-3))=\frac{(2\alpha^2-\alpha)n-6 \alpha^2+ \alpha+1}{ \alpha}\geq0
$$
for $n \geq 8$ and $\alpha \in[\frac{1}{2},1)$, then we have
$\rho_\alpha(G)\leq\alpha(n-2)+\frac{(1-\alpha)^2}{\alpha}.$
Combining this with (\ref{eq4}), we have
$
\rho_\alpha(G)<\rho_\alpha\left(K_{2, n-2}\right)
$
for $n \geq 8$ and $\alpha \in[\frac{1}{2},1)$.

\textbf{Subcase 1.2.} $e(N(w))=1$.

By Lemma \ref{chord}, we know that each vertex in $B$ exactly has one neighbor in $N(w)$. It follows that $e(N(w),B)<n-3$. Thus, $\sum_{wv \in E(G)} d(v)= d(w) + 2e(N(w)) + e(N(w),B)<n+1$. Combining this with (\ref{eq5}), we have
$$
\rho_\alpha(G) \leq 2 \alpha+\frac{1-\alpha}{2} \sum_{w v \in E(G)} d(v) < 2 \alpha+\frac{(1-\alpha)(n+1)}{2}.
$$
Noting that
$$
\alpha(n-2)+\frac{(1-\alpha)^2}{\alpha}-(2 \alpha+\frac{(1-\alpha)(n+1)}{2})=\frac{(3\alpha^2-\alpha)n-5\alpha^2-5 \alpha+2}{ \alpha}\geq0
$$
for $n \geq 8$ and $\alpha \in[\frac{1}{2},1)$, then we have
$\rho_\alpha(G)<\alpha(n-2)+\frac{(1-\alpha)^2}{\alpha}.$
Combining this with (\ref{eq4}), we have
$
\rho_\alpha(G)<\rho_\alpha\left(K_{2, n-2}\right)
$
for $n \geq 8$ and $\alpha \in[\frac{1}{2},1)$.

%%case2
\textbf{Case 2.} $3\leq d(w) \leq n-3$.

In order to prove $\rho_\alpha(G)<\rho_\alpha(K_{2,n-2})= \frac{1}{2}(\alpha n+\sqrt{\alpha^2 n^2+8(n-2)(1-2\alpha)})$,
it is enough to prove $\rho_\alpha(G)^2-\alpha n \rho_\alpha(G) +2(2\alpha-1)(n-2)< 0.$
For convenience, we denote $A_\alpha(G)=A_\alpha$, $A(G)=A$ and $D(G)=D$. Let
$$
B=(b_{ij})_{n\times n}=A_\alpha^2-\alpha n A_\alpha +2(2\alpha-1)(n-2)I_n,
$$
where $I_n$ is the $n\times n$ identity matrix. Let $c_u(B)$ be the sum of all elements in the $u$-th column of $B$. Then we have the following claim.
\par{\textbf{Claim 2.1.} $c_u(B)< 0$ for $n\geq 8$ and $\alpha \in[\frac{1}{2},1)$.
\par{\textbf{Proof.}
Since $A_\alpha=\alpha D+(1-\alpha)A$, then
$$
\begin{aligned}
B=&(\alpha D+(1-\alpha) A)^2-\alpha n(\alpha D+(1-\alpha) A)+2(2\alpha-1)(n-2) I_n \\
=& \alpha^2 D^2+(1-\alpha)^2 A^2+\alpha(1-\alpha) D A+\alpha(1-\alpha) A D-\alpha^2 n D
-(\alpha n-\alpha^2 n) A\\&+2(2\alpha-1)(n-2)I_n.
\end{aligned}
$$

It is easy to see that $c_u(A)=c_u(D)=d(u)$, $c_u(A^2)=c_u(DA)=\sum_{u v \in E(G)} d(v)$ and $c_u(AD)=d^2(u)$. By Lemma \ref{mn}, we have
$\sum_{w v \in E(G)} d(v)< |E(G)|\leq 2n-2$. It follows that
$$
\begin{aligned}
c_u(B)=& \alpha^2 d^2(u)+(1-\alpha)^2 \sum_{uv \in E(G)} d(v)+\alpha(1-\alpha) \sum_{uv \in E(G)} d(v)+\alpha(1-\alpha)d^2(u)-\alpha^2 n d(u)
\\&-(\alpha n-\alpha^2 n) d(u)+2(2\alpha-1)(n-2) \\
=&\alpha d^2(u)+(1-\alpha) \sum_{uv \in E(G)} d(v)-\alpha n d(u)+2(2\alpha-1)(n-2)\\
< & \alpha d^2(u)+(1-\alpha) (2n-2)-\alpha n d(u)+2(2\alpha-1)(n-2) \\
=& \alpha(d^2(u)-n d(u)+2 n-6)+2 \\
\leq & \max \left\{\alpha(9-3n+2 n-6), \alpha((n-3)^2-n(n-3)+2 n-6)\right\}+2 \\
=& \alpha(-n+3)+2< 0
\end{aligned}
$$
for $n\geq 8$ and $\alpha \in[\frac{1}{2},1)$.

This completes the proof of the claim. $\qedsymbol$}
}
\par{Let $X=\left(x_1, x_2, \ldots, x_n\right)^{T}$ be the $\alpha$-Perron vector of $G$ corresponding to $\rho_\alpha(G)$ satisfying $\sum_{i=1}^n x_i=1$. Then
$$
B X=(\rho_\alpha(G)^2-\alpha n \rho_\alpha(G) +2(2\alpha-1)(n-2))X.
$$
Hence we have
$$
\begin{aligned}
&\rho_\alpha(G)^2-\alpha n \rho_\alpha(G) +2(2\alpha-1)(n-2)\\
&=\sum_{i=1}^n(\rho_\alpha(G)^2-\alpha n \rho_\alpha(G) +2(2\alpha-1)(n-2)) x_i \\
&=\sum_{i=1}^n(B X)_i=\sum_{i=1}^n(\sum_{j=1}^n b_{i j} x_j) =\sum_{j=1}^n(\sum_{i=1}^n b_{i j}) x_j=\sum_{j=1}^n c_j(B) x_j<0.
\end{aligned}
$$}
\par{Combining the above arguments, we have $\rho_\alpha(G)\leq \rho_\alpha(K_{2,n-2})$ for $n\geq 8$ and $\alpha \in[\frac{1}{2},1)$, with equality if and only if $G\cong K_{2,n-2}$.

These complete the proof. $\qedsymbol$}
%%% End proof Theorem 1.1 %%%
%%%%%%%%%%%%%%%%%%%%%%%%%%%%%%%%%%%%%%%%%%%%%%%%%%%%%%%%%%%%%%%%%%%%%%%%%%%%%%

%%%%%%%%%%%%%%%%%%%%%%%%%%%%%%%%%%%%%%%%%%%%%%%%%%%%%%%%%%%%%%%%%%%%%%%%%%%%%%
%%% Part 4: Proof Theorem 1.2 %%%
\section{Proof of Theorem \ref{th1.2}}
\label{sec:ch-sufficient}
\noindent In this section, we give the proof of Theorem \ref{th1.2}.\\ \hspace*{\fill} \\
\textbf{Proof.}
Let $G$ be a minimally 3-connected graph with order $n$.
Notice that $W_n$ is a minimally 3-connected graph. By Lemma \ref{lb}, we have
\begin{equation}{\label{eq6}}
\rho_{\alpha}(W_n)> \alpha\Delta(W_n)+\frac{(1-\alpha)^{2}}{\alpha}
=\alpha(n-1)+\frac{(1-\alpha)^{2}}{\alpha}
\end{equation}
for $\alpha\in[\frac{1}{2},1)$.
Let $w$ be a vertex of $G$ such that
$$\alpha d(w)+\frac{1-\alpha}{d(w)} \sum_{wv \in E(G)}d(v)=\max_{u \in V(G)}\left\{\alpha d(u)+\frac{1-\alpha}{d(u)} \sum_{uv \in E(G)} d(v)\right\}.$$\\
By Lemma \ref{ub1}, we have
\begin{equation}{\label{eq7}}
\rho_{\alpha}(G)\leq \alpha d(w)+\frac{1-\alpha}{d(w)} \sum_{wv \in E(G)}d(v).
\end{equation}

Denote by $B=V(G)\backslash N[w]$. We assert that $d_{N(w)}(u)\leq2$ for $u\in N(w)$. Otherwise, we assume that $d_{N(w)}(u_1)\geq3$ for $u_1\in N(w)$
and $u'_{1}\in N_{N(w)}(u_1)$. Then there exists a cycle $wu_1u'_{1}w$ with $d(w)\geq4$ and $d(u_1)\geq4$, which contradicts Lemma \ref{degree3}. Hence, we have $d_{N(w)}(u)\leq2$ for $u\in N(w)$.
If $\Delta(G)=n-1$, then $G[N(w)]$ is the union of some cycles by Lemma \ref{degree}. It follows that $G\cong W_n$. Otherwise $w$ is a cut vertex. If $\Delta(G)\leq n-2$, then $3\leq d(w)\leq n-2$ by Lemma \ref{degree}. We have the following claim.\\
\textbf{Claim 1.} $G[N(w)]$ contains no cycles for $3\leq d(w)\leq n-2$. \\
\textbf{Proof.}
If there exists a cycle $C_q$ in $G[N(w)]$, then $W_{q+1}\subseteq G$. Since $3\leq d(w)\leq n-2$, we have $3\leq q\leq n-2$, then $|V(W_{q+1})|< n$. Notice that every vertex of $W_{q+1}$ is contained in some $K_3$. By Lemma \ref{degree3},
there exists exactly one vertex $u\in V(W_{q+1})$ such that $d(u)\geq4$.
It follows that $u$ is a cut vertex, which contradicts the fact that $G$ is a minimally 3-connected graph. $\qedsymbol$

By Claim 1, we know that $G[N(w)]$ is the union of some paths and isolated vertices. We consider the following two cases.

\textbf{Case 1.} $e(N(w))=0$.

In this case, we have
$
e(N(w),B)\leq d(w) |B|=d(w)(n-d(w)-1).
$ It follows that
$$
\sum_{w v \in E(G)} d(v)=2e(N(w))+d(w)+e(N(w),B)\leq d(w)(n-d(w)).
$$
Combining with (\ref{eq7}), we have
$$
\begin{aligned}
\rho_{\alpha}(G)&\leq \alpha d(w)+(1-\alpha)(n-d(w))\\
&=(2\alpha-1)d(w)+(1-\alpha)n\leq(n-4)\alpha+2.
\end{aligned}
$$
Noting that
$$
\alpha(n-1)+\frac{(1-\alpha)^2}{\alpha}-((n-4)\alpha+2)=\frac{(2\alpha-1)^2}{\alpha}\geq0
$$ for $\alpha \in[\frac{1}{2},1)$, we have $\rho_\alpha(G)\leq \alpha(n-1)+\frac{(1-\alpha)^2}{\alpha}$. Combining this with (\ref{eq6}), we have
$$
\rho_\alpha(G)\leq \alpha(n-1)+\frac{(1-\alpha)^2}{\alpha}<\rho_\alpha(W_n)
$$for $\alpha\in[\frac{1}{2},1)$.

\textbf{Case 2.} $e(N(w))\geq1$.

We consider the following two subcases.

\textbf{Subcase 2.1.} $d(w)=3$.

In this case, $|B|=n-4$. If $e(N(w))=1$, by Lemma \ref{degree3}, we have $e(N(w),B)\leq(d(w)-1)|B|+1$. If $e(N(w))=2$, we assume $P_3=u_1u_2u_3\in G[N(w)]$, then $N_B(u_2)\neq \varnothing$.
Otherwise, $G-\{u_1, u_3\}$ is disconnected, which contradicts the fact that $G$ is a minimally 3-connected graph. By Lemma \ref{degree3}, we have $e(N(w),B)\leq 2+(d(w)-2) |B|$. Hence we have $e(N(w),B)\leq(3-e(N(w))) |B|+e(N(w))$ for $e(N(w))=1,2$. It follows that
$$
\begin{aligned}
\sum_{w v \in E(G)} d(v)&=2e(N(w))+d(w)+e(N(w),B)\\&\leq 2e(N(w))+3+(3-e(N(w)))(n-4)+e(N(w))\\
&=(3-e(N(w)))n+7e(N(w))-9.
\end{aligned}
$$
Combining with (\ref{eq7}), we have
$$
\begin{aligned}
\rho_{\alpha}(G)&\leq 3\alpha +\frac{1-\alpha}{3}((3-e(N(w)))n+7e(N(w))-9)\\
&\leq max\left\{3\alpha +\frac{1-\alpha}{3}(2n-2),3\alpha +\frac{1-\alpha}{3}(n+5)\right\}\\
&=3\alpha +\frac{2}{3}(1-\alpha)(n-1).
\end{aligned}
$$
Noting that
$$
\alpha(n-1)+\frac{(1-\alpha)^2}{\alpha}-(3\alpha +\frac{2}{3}(1-\alpha)(n-1))=\frac{(5\alpha^2-2\alpha)n-11\alpha^2-4\alpha+3}{3\alpha}\geq0
$$ for $n\geq7$ and $\alpha \in[\frac{1}{2},1)$, we have $\rho_\alpha(G)\leq \alpha(n-1)+\frac{(1-\alpha)^2}{\alpha}$. Combining this with (\ref{eq6}), we have
$$
\rho_\alpha(G)\leq \alpha(n-1)+\frac{(1-\alpha)^2}{\alpha}<\rho_\alpha(W_n)
$$for $n\geq7$ and $\alpha\in[\frac{1}{2},1)$.

\textbf{Subcase 2.2.} $4\leq d(w)\leq n-2$.

Let $P_{t_i}\;(1\leq i\leq s,\;t_i\geq2)$ be the s paths in $G[N(w)]$. Then $e(N(w))=\sum^{s}_{i=1}(t_i-1)$. We assert that $|B|\geq2$. Otherwise, if $|B|=1$, we assume that $B=\{v\}$.
Denote by $c$ the number of components in $G[N(w)]$. If $c=1$, then $P_{d(w)}$ is the only one path in $G[N(w)]$. By Lemma \ref{degree3}, we have $d(v)=2$, which contradicts Lemma \ref{degree}.
If $c\geq2$, then $G-\{w,v\}$ is disconnected, which contradicts the fact that $G$ is a minimally 3-connected graph. Hence, we have $|B|\geq2$.
By Lemma \ref{degree3}, we have $e(N(w),B)\leq2s+(d(w)-\sum^{s}_{i=1}t_i)|B|$. Hence we have
$$
\begin{aligned}
\sum_{w v \in E(G)} d(v)&=2e(N(w))+d(w)+e(N(w),B)\\
&\leq 2\sum^{s}_{i=1}(t_i-1)+d(w)+2s+(d(w)-\sum^{s}_{i=1}t_i)|B|\\
&=d(w)(|B|+1)+(2-|B|)\sum^{s}_{i=1}t_i.
\end{aligned}
$$
Combining with (\ref{eq7}) and $|B|\geq2$, we have
$$
\begin{aligned}
\rho_{\alpha}(G)&\leq \alpha d(w) +\frac{1-\alpha}{d(w)}(d(w)(|B|+1)+(2-|B|)\sum^{s}_{i=1}t_i)\\
&\leq\alpha d(w)+(1-\alpha)(|B|+1)\\
&=\alpha d(w)+(1-\alpha)(n-d(w))\\
&\leq(n-4)\alpha+2.
\end{aligned}
$$
Noting that
$$
\alpha(n-1)+\frac{(1-\alpha)^2}{\alpha}-((n-4)\alpha+2)=\frac{(2\alpha-1)^2}{\alpha}\geq0
$$ for $\alpha \in[\frac{1}{2},1)$, we have $\rho_\alpha(G)\leq \alpha(n-1)+\frac{(1-\alpha)^2}{\alpha}$. Combining this with (\ref{eq6}), we have
$$
\rho_\alpha(G)\leq \alpha(n-1)+\frac{(1-\alpha)^2}{\alpha}<\rho_\alpha(W_n)
$$for $\alpha\in[\frac{1}{2},1)$.
\par{Combining the above arguments, we have $\rho_\alpha(G)\leq \rho_\alpha(W_n)$ for $n\geq 7$ and $\alpha \in[\frac{1}{2},1)$, with equality if and only if $G\cong W_n$.}

These complete the proof. $\qedsymbol$

%%% Part 4 end %%%
%%%%%%%%%%%%%%%%%%%%%%%%%%%%%%%%%%%%%%%%%%%%%%%%%%%%%%%%%%%%%%%%%%%%%%%%%%%%%%
\section*{Data availability}
\par{No data was used for the research described in the article.}
\section*{Conflict of interest}
\par{The authors declare that they have no conflict of interest.}
%\section*{Acknowledgments}
%\par{The authors sincerely thank the anonymous referees for their valuable comments and constructive suggestions on the original manuscript.}
%%%%%%%%%%%%%%%%%%%%%%%%%%%%%%%%%%%%%%%%%%%%%%%%%%%%%%%%%%%%%%%%%%%%%%%%%%%%%%

%%% Part 5 end %%%%
%%%%%%%%%%%%%%%%%%%%%%%%%%%%%%%%%%%%%%%%%%%%%%%%%%%%%%%%%%%%%%%%%%%%%%%%%%%%%%

\end{document}